\documentclass[12pt]{amsart}
\usepackage{amsmath,amsthm,amsfonts,amssymb}
\begin{document} 
\newcommand{\B}{{\mathbb B}}
\newcommand{\C}{{\mathbb C}}
\newcommand{\N}{{\mathbb N}}
\newcommand{\Q}{{\mathbb Q}}
\newcommand{\Z}{{\mathbb Z}}
\renewcommand{\P}{{\mathbb P}}
\newcommand{\R}{{\mathbb R}}
\newcommand{\rc}{\subset}
\newcommand{\rank}{\mathop{rank}}
\newcommand{\trace}{\mathop{tr}}
\newcommand{\dimc}{\mathop{dim}_{\C}}
\newcommand{\Lie}{\mathop{Lie}}
\newcommand{\Auto}{\mathop{{\rm Aut}_{\mathcal O}}}
\newcommand{\Aut}{\mathop{{\rm Aut}}}
\newcommand{\alg}[1]{{\mathbf #1}}

\newtheorem*{definition}{Definition}
\newtheorem*{claim}{Claim}
\newtheorem{corollary}{Corollary}
\newtheorem*{Conjecture}{Conjecture}
\newtheorem*{SpecAss}{Special Assumptions}
\newtheorem{example}{Example}
\newtheorem*{remark}{Remark}
\newtheorem*{observation}{Observation}
\newtheorem*{fact}{Fact}
\newtheorem*{remarks}{Remarks}
\newtheorem{lemma}{Lemma}
\newtheorem{proposition}{Proposition}
\newtheorem{theorem}{Theorem}


\title{Holomorphic discs with dense images}

%
%
\author{Franc Forstneri\v c}
\address{Institut of mathematics, physics and mechanics, 
University of Ljubljana, Jadranska 19, 1000 Ljubljana, Slovenia}
\email{franc.forstneric@fmf.uni-lj.si}
\thanks{Research of the first author was supported by grants P1-0291
and J1-6173, Republic of Slovenia.}

%
%
\author{J\"org Winkelmann}
\address{Institut Elie Cartan (Math\'ematiques), Universit\'e 
Henri Poincar\'e Nancy 1, B.P.\ 239, F-54506 Vand\oe uvre-les-Nancy Cedex, France}
\email{jwinkel@member.ams.org\newline\indent{\itshape Webpage: }%
http://www.math.unibas.ch/\~{ }winkel/
}

%
%
\date{\today}

%
%
\begin{abstract}
Let $\Delta$ be the open unit disc in $\C$,
$X$ a connected complex manifold and $\mathcal D$ the set of all
holomorphic maps $f \colon \Delta\to X$ with $\overline{f(\Delta)}=X$.
We prove that $\mathcal D$ is dense in  $Hol(\Delta,X)$.
\end{abstract}

\maketitle

%
%
%
%
\section{Introduction}
Let $\Delta_r=\{z\in\C\colon |z|<r\}$ and $\Delta=\Delta_1$.
In \cite{W} the second author proved that 
{\em for any irreducible complex space $X$ there exists a
holomorphic map $\Delta \to X$ with dense image}, and 
he raised the question whether the set of all holomorphic maps
$\Delta\to X$ with dense image 
forms a dense subset of the set $Hol(\Delta,X)$ of all holomorphic maps 
$\Delta\to X$ with respect to the topology of locally uniform convergence.

In this paper we show that the answer to this question is positive
if $X$ is smooth, but negative for some singular space.

\begin{theorem}
\label{t1}
For any connected complex manifold $X$ the set of holomorphic maps 
$\Delta\to X$ with dense images forms a dense subset 
in $Hol(\Delta,X)$. The conclusion fails for some
singular complex surface $X$.
\end{theorem}

The situation is quite different for {\em proper discs},
i.e., proper holomorphic maps $\Delta\to X$.
The paper \cite{FG} contains an example of a 
non-pseudoconvex bounded domain $X\subset \C^2$
such that a certain nonempty open subset $U\subset X$
is not intersected by the image of any proper 
holomorphic disc $\Delta\to X$.
On the other hand, proper holomorphic discs exist in great 
abundance in {\em Stein manifolds} \cite{G}, \cite{BD1}, \cite{BD2}.

%
%
%
%
\section{Preparations}

\begin{lemma}\label{lemma1}
Let $W_n$ be a decreasing sequence (i.e., $W_{n+1}\subset W_n$)
of open sets with $\Delta\subset W_n\subset\Delta_2$
for every $n$. Let $K=\cap_n \overline{W}_n$ and assume that 
the interior of $K$ coincides with $\Delta$. Furthermore assume that 
there are biholomorphic maps $\phi_n\colon \Delta\to W_n$
with $\phi_n(0)=0$ for $n=1,2,\ldots$.

Then there exists an automorphism $\alpha\in\Aut(\Delta)$
and a subsequence $(\phi_{n_k})$ of the sequence
$(\phi_n)$ such that $\phi_{n_k}\circ\alpha^{-1}$ converges locally
uniformly to the identity map $id_\Delta$ on $\Delta$.
\end{lemma}

\begin{proof}
Montel's theorem shows that, after passing to a suitable subsequence,
we have $\lim_{n\to\infty} \phi_n=\alpha\colon \Delta\to K$ and
$\lim_{n\to\infty} (\phi_n^{-1}|_\Delta)=\beta\colon \Delta\to \overline\Delta$.
Since the limit maps are holomorphic and satisfy
$\alpha(0)=0$ and $\beta(0)=0$, we conclude that 
$\alpha(\Delta)\subset {\rm Int}K=\Delta$ and 
$\beta(\Delta)\subset \Delta$. Moreover 
$\alpha\circ\beta=id_\Delta=\beta\circ\alpha$,
and hence both $\alpha$ and $\beta$ are automorphisms 
of $\Delta$ (indeed, rotations $z\to z e^{it}$).
\end{proof}

We also need the following special case of 
a result of the first author (theorem 3.2 in \cite{F}):

\begin{proposition}
\label{propf}
Let $X$ be a complex manifold, $0<r<1$, 
$E$ the real line segment $[1,2]\subset \C$,
$K=\overline \Delta\cup E$, $U$ an open neighbourhood 
of $\overline \Delta$ in $\C$, $S$ a finite subset of $K$ and 
$f\colon U\cup E \to X$ a continuous map which 
is holomorphic on $U$.

Then there is a sequence of pair of open neighbourhoods
$W_n \subset \C$ of $K$ and holomorphic maps $g_n\colon W_n\to X$
such that:
\begin{enumerate}
\item $g_n|_K$ converges uniformly to $f|_K$ as $n\to\infty$, and 
\item $g_n(a)=f(a)$ for all $a\in S$ and $n\in \N$.
\end{enumerate}
\end{proposition}

\section{Towards the main result}
In this section we prove the following proposition
which is the main technical result in the paper.
The first statement in theorem \ref{t1} (\S 1) is an 
immediate corollary.

\begin{proposition}
\label{mainprop}
Let $X$ be a connected complex manifold endowed with a complete
Riemannian metric and induced distance $d$, 
$S$ a countable subset of $X$, $f\colon \Delta\to X$ 
a holomorphic map, $\epsilon>0$ and $0<r<1$.

Then there exists a holomorphic map $F\colon \Delta\to X$
such that
\begin{enumerate}
\item[(a)] $S\subset F(\Delta)$, and
\item[(b)] $d(f(z),F(z))\le\epsilon$ for all $z\in\Delta_r$.
\end{enumerate}
\end{proposition}

\begin{proof}
Let $s_1,s_2,s_3,..$ be an enumeration of the elements of $S$.
We shall inductively construct a sequence of 
holomorphic maps $f_n\colon \Delta\to X$,
numbers $r_n\in (0,1)$ and points
$a_{1,n},\ldots,a_{n,n} \in \Delta$
satisfying the following properties for $n=0,1,2,\ldots$:
\begin{enumerate}
\item
$f_0=f$ and $r_0=r$,
\item
$(r_n+1)/2<r_{n+1}<1$,
\item
$f_n(a_{j,n})=s_j$ for $n\ge 1$ and $j=1,2,\ldots,n$,
\item
$d(f_n(z),f_{n+1}(z))<2^{-(n+1)}\epsilon$ for all $z\in\Delta_{r_n}$, and
\item
$d_\Delta (a_{j,n},a_{j,n+1})<2^{-n}$ for $j=1,2,\ldots,n$ where 
$d_\Delta$ denotes the Poincar\'e distance on the unit disc.
\end{enumerate}

Assume inductively that the data for level $n$ 
(i.e., $f_n$, $r_n$, $a_{j,n}$) have been chosen. 
(For $n=0$ we do not have any points $a_{j,0}$.)
With $n$ fixed we choose an increasing sequence of real numbers 
$\lambda_k$ with $\lambda_k>r_n$ and $\lim_{k\to\infty}\lambda_k=1$.
For every $k\in \N$ the map 
$\widetilde g_k(z) \stackrel{def}{=} f_n(\lambda_k z) \in X$
is defined and holomorphic on the disc 
$\Delta_{1/\lambda_k} \supset \overline \Delta$. 
After a slight shrinking of its domain we can 
extend it continuously to the segment
$E=[1,2] \subset \C$ such that the right end point 
$2$ of $E$ is mapped to the next point $s_{n+1}\in S$ 
(this is possible since $X$ is connected). 

Applying proposition \ref{propf} to the extended map 
$\widetilde g_k$ we obtain for every $k\in \N$ 
an open neighbourhood 
$V_k\subset\C$ of $K=\overline \Delta \cup E$ 
and a holomorphic map $g_k\colon V_k\to X$ such that
\begin{enumerate}
\item[(i)]
$|g_k(z)-f_n(\lambda_kz)|<2^{-k}$ for all $z\in\overline\Delta$, 
\item[(ii)]
$g_k(2)=s_{n+1}$, and 
\item[(iii)]
$g_k(a_{j,n}/\lambda_k)=f_n(a_{j,n})=s_j$ 
for $j=1,\ldots,n$.
\end{enumerate}

Next we choose a decreasing sequence of simply connected open sets 
$W_k \subset\C$ $(k\in\N)$ with $K \subset W_k\subset V_k$ and $K=\cap_k \overline W_k$.
Notice that ${\rm Int} K=\Delta$.
By lemma~\ref{lemma1} there is a sequence of biholomorphic maps 
$\phi_k\colon \Delta\to W_k$ with $\lim_{k\to\infty} \phi_k=id_\Delta$.

Consider the holomorphic maps $h_k=g_k\circ\phi_k\colon \Delta\to X$.
By our construction we know that $\lim_{k\to\infty}h_k=f_n$ 
locally uniformly on $\Delta$.

To fulfill the inductive step it thus suffices to 
choose $f_{n+1}=h_k$ for a sufficiently large $k$,
$a_{j,n+1}=a_{j,n}/\lambda_k$ ($j=1,\ldots,n$), 
$a_{n+1,n+1}=\phi_k^{-1}(2)$. Finally we choose a number 
$r_{n+1}$ satisfying 
\[
\max\{|a_{n+1,n+1}|, \frac{r_n+1}{2}\} < r_{n+1} <1.
\]
This completes the inductive step.

By properties (2) and (4) the sequence $f_n$ converges 
locally uniformly in $\Delta$ to a holomorphic
map $F\colon \Delta\to X$. Aided by property (1) we 
also control $d(f(z),F(z))$ for $z\in\Delta_r$.
Since the Poincar\'e metric is complete, 
property (5) insures that for every 
fixed $j\in \N$ the sequence $a_{j,n} \in\Delta$ 
$(n=j,j+1,\ldots)$ has an accumulation point $b_j$ 
inside of $\Delta$, and (3) implies $F(b_j)=s_j$ 
for $j=1,2,\ldots$. Hence $S\subset F(\Delta)$. 
\end{proof}

\section{Singular spaces}
We use an example of Kaliman and Zaidenberg \cite{KZ} to show 
that for a complex spaces $X$ with singularities the set of maps
$\Delta\to X$ with dense image need not be dense
in $Hol(\Delta,X)$. We denote by $Sing(X)$ the singular 
locus of $X$.

\begin{proposition}
There is a singular compact complex surface $S$, 
a non-constant holomorphic map 
$f\colon \Delta\to S$ and an open neighbourhood $\Omega$ of $f$
in $Hol(\Delta,S)$ such that $g(\Delta)\subset Sing(S)$
for every $g\in\Omega$.
\end{proposition}

\begin{proof}
In \cite{KZ} Kaliman and Zaidenberg constructed an example of a singular surface
$S$ with normalization $\pi\colon Z\to S$ such that $S$ contains a 
rational curve $C\simeq \P^1$ while $Z$ is smooth and hyperbolic.
Denote by $d_Z$ the Kobayashi distance function on $Z$.
We choose two distinct points $p,q\in C$ and open 
relatively compact neighbourhoods $V$ of $p$ 
and $W$ of $q$ in $S$ such that 
$\overline V\cap\overline W=\emptyset$.
The preimages $\pi^{-1}(\overline{V})$ and
$\pi^{-1}(\overline{W})$ in $Z$ are also compact,
and since $Z$ is hyperbolic we have 
\[
r=\min\{ d_Z(x,y)\colon x\in \pi^{-1}(\overline{V}), y\in \pi^{-1}(\overline{W} )\} > 0.
\]

Fix a point $a\in\Delta$ with $0<d_\Delta(0,a)<r$ and let
$\Omega$ consist of all holomorphic maps $g\colon \Delta\to S$
satisfying $g(0)\in V$ and $g(a)\in W$. Since both $p$
and $q$ are lying on the rational curve $C$, 
there is a holomorphic map $g\colon\Delta\to C$
with $g(0)=p\in V$ and $g(a)=q\in W$; hence the set $\Omega$
is not empty. Clearly $\Omega$ is open in $Hol(\Delta,S)$.

To conclude the proof it remains to show that 
$g(\Delta) \subset Sing(S)$ for all $g\in\Omega$.
Indeed, a holomorphic map $g\colon \Delta \to S$ 
with $g(\Delta)\not\subset Sing(S)$ admits 
a holomorphic lifting $\widetilde g\colon \Delta\to Z$
with $\pi\circ \widetilde g=g$. If $g\in \Omega$ then by 
construction 
\[
 d_Z(\widetilde g(0),\widetilde g(a)) \ge r > d_\Delta(0,a)
\] 
which violates the distance decreasing property for 
the Kobayashi pseudometric. This contradiction establishes the claim.
\end{proof}

In particular, we see that in this example the set of all holomorphic
maps $f\colon \Delta\to S$ with dense image {\em does not} constitute
a dense subset of $Hol(\Omega,S)$.

\end{document}